\documentclass[twocolumn]{autart}    % Enable this line and disable the 
\usepackage{cite}
\usepackage{url}
\usepackage{amssymb,amsfonts}
\usepackage{multirow}
\usepackage{textcomp}
\usepackage{epsfig}

\usepackage{epstopdf}
\usepackage[ruled]{algorithm}
\usepackage{enumerate}
\usepackage{float}
\usepackage{graphics,graphicx}
\usepackage{subfigure}
\usepackage{algpseudocode}

\usepackage{array}
\usepackage{mathtools}
\usepackage[justification=centering]{caption}
\allowdisplaybreaks[3]
\usepackage{url}
\DeclareGraphicsExtensions{.png}
\graphicspath{{./Pics/}}
\usepackage{color}
\usepackage{xspace}

\usepackage{hyperref}
\hypersetup{
    colorlinks=true,
    linkcolor=blue,
    filecolor=gray,      
    urlcolor=blue,
    citecolor=blue,
}

\def\BibTeX{{\rm B\kern-.05em{\sc i\kern-.025em b}\kern-.08em
		T\kern-.1667em\lower.7ex\hbox{E}\kern-.125emX}}

\newtheorem{assumption}{Assumption}

\begin{document}

\begin{frontmatter}

\title{A Distributed Gradient-based Algorithm for Optimization Problems with Coupled Equality Constraints \thanksref{mytitlenote}}
\thanks[mytitlenote]{This work was supported in part by the National Science Foundation under the grant number CMMI–2243930. The material in this paper was
not presented at any conference. Corresponding author: Zongli Lin.}
%\thanksref{footnoteinfo}} % Title, preferably not more 

\author[address]{Chenyang Qiu}\ead{nzp4an@virginia.edu},
\author[address]{Zongli Lin}\ead{zl5y@virginia.edu}

\address[address]{Charles L. Brown Department of Electrical and Computer Engineering, University of Virginia, Charlottesville, VA 22904, USA}

% \author[Paestum]{Chenyang Qiu}\ead{nzp4an@virginia.edu},    % Add the 
% % \author[Baiae]{Chenyang Qiu}\ead{nzp4an@virginia.edu}  % (ead) as shown

% \address[Paestum]{the Charles L. Brown Department of Electrical and Computer Engineering, University of Virginia, Charlottesville, VA 22904, USA}  % Please supply                                              
% \address[Rome]{Senate House, Rome}             % full addresses
% \address[Baiae]{The White House, Baiae}        % here.


\begin{abstract}
This paper studies a class of distributed optimization problems with coupled equality constraints in networked systems. Many existing distributed algorithms rely on solving local subproblems via the $\operatorname{argmin}$ operator in each iteration. Such approaches become computationally burdensome or intractable when local cost functions are complex. To address this challenge, we propose a novel distributed gradient-based algorithm that avoids solving a local optimization problem at each iteration by leveraging first-order approximations and projection onto local feasible sets. The algorithm operates in a fully distributed manner, requiring only local communication without exchanging gradients or primal variables. We rigorously establish sublinear convergence for general convex cost functions and linear convergence under strong convexity and smoothness conditions. Numerical simulation on the IEEE 118-bus system demonstrates the superior computational efficiency and scalability of the proposed method compared to several state-of-the-art distributed optimization algorithms.

\end{abstract}

\end{frontmatter}

\section{Introduction}\label{sec:introduction}

% \cite{xu2018dual}
The optimization problem with coupled equality constraints, also known as the resource allocation or economic dispatch problem, serves as a fundamental model in various real-world applications, including efficient power system operation \cite{xie2024virtual}, channel assignment in wireless communication networks \cite{zhao2018simple}, and task scheduling \cite{jamil2022resource}.
At its core, the optimization problems with coupled equality constraints seek to optimally allocate available resources to satisfy demand at the minimum possible cost, while adhering to a range of operational constraints.
% Such problems are formulated as a constrained optimization problem that minimizes the total generation cost subject to demand satisfaction and system limitations \cite{yang2016distributed}.
Historically, such problems have been solved via centralized approaches, where a central coordinator, equipped with global system information, computes the optimal solution \cite{hetzer2008economic}. While effective, these centralized methods incur significant communication overhead, raise privacy concerns, and are vulnerable to the failure of the central node. In the meantime, the rapid development of several advanced technologies, such as distributed energy resources, cyber-networks, and cloud computing, is making conventional centralized management strategies increasingly inadequate \cite{akorede2010distributed, muhtadi2021distributed}.

These practical limitations of centralized methods, together with the growing deployment of advanced technologies, have spurred the development of distributed algorithms.
In the distributed framework, each agent conducts local optimization and communicates with neighboring units through limited information exchange. This decentralized strategy significantly improves scalability, enhances system resilience, and offers better privacy preservation. Consensus-based algorithms have been extensively adopted to solve the optimization problem with coupled equality constraints in a distributed manner. These algorithms exploit local communication to iteratively achieve consensus on marginal costs, enabling fully distributed and scalable solutions without the need for centralized coordination \cite{wang2018distributed, liu2020constrained, chen2024privacy, ji2024self}. However, existing consensus-based approaches often assume quadratic cost functions and relatively simple constraints, which limit their applicability to more general settings involving non-quadratic objectives.

To handle more general or nonsmooth cost functions beyond quadratic functions, a range of advanced distributed optimization algorithms has been proposed. A common approach involves reformulating the distributed optimization problem with coupled equality constraints by introducing two convex sets with associated indicator functions, thereby enabling variable splitting and problem decomposition. Building on the Alternating Direction Method of Multipliers (ADMM) framework, several works have made notable progress \cite{chang2014multi, chen2017admm, sun2024distributed}. Ref. \cite{chang2014multi} provided a class of algorithms that solve the distributed optimization problem with coupled equality constraints with convergence guarantees. Further, Ref. \cite{chen2017admm} introduced a fully distributed ADMM-based algorithm with convergence guarantees.
Ref. \cite{su2021distributed} proposes a distributed primal–dual algorithm, where each iteration embeds a finite-time multi-step consensus iterations to get local estimates of the global constraint violation.
More recently, Ref. \cite{falsone2023augmented} proposed an augmented Lagrangian tracking method to address distributed optimization problems involving coupled constraints, broadening the applicability of ADMM-based strategies.

% To address privacy concerns in communication, Ref. \cite{sun2024distributed} developed a privacy-enhanced ADMM algorithm for microgrid economic dispatch, which integrates encryption–decryption protocols into the ADMM framework to ensure secure communication without compromising convergence performance.

Another line of research tries to solve the distributed optimization problem with coupled equality constraints from a duality perspective, by reformulating it as a consensus optimization problem and solving its dual formulation \cite{nedic2018improved, xu2018dual, zhang2020distributed, wu2025distributed, li2019distributed}. For instance, based on the mirror relationship between the consensus optimization problem and the optimization problem with coupled equality constraints, Ref. \cite{nedic2018improved} proposed three algorithms, Mirror-P-EXTRA, Mirror-EXTRA, and Mirror-PG-EXTRA. Ref. \cite{xu2018dual} utilized the first-order information of the objective function and provided convergence rates for both convex and strongly convex cases. Leveraging the duality, Ref. \cite{zhang2020distributed} developed a distributed dual gradient tracking algorithm for resource allocation. Under the assumption of strong convexity and in the absence of local constraints, the algorithm achieves ergodic sublinear convergence to the Karush–Kuhn–Tucker (KKT) condition at a rate of $O(1/k)$. Algorithms proposed in Ref. \cite{wu2025distributed} achieve an ergodic convergence rate of $O(1/k)$ considering local constraints and have been extended for differential privacy. 

It is worth noting that many existing distributed algorithms for solving distributed optimization problems with coupled equality constraints rely heavily on the use of the $\operatorname{argmin}$ operator and require each agent to solve a local optimization problem at every iteration. While this approach can be efficient when local objective functions are simple, such as in the case of quadratic functions with closed-form solutions, its applicability becomes limited when dealing with more complex functions, such as exponential type cost functions \cite{he2019optimizing}, where solving the local optimization problem at each iteration is computationally expensive or even intractable.

To address this challenge, we propose a novel distributed algorithm, termed the Distributed Gradient-based Algorithm (DGA), that circumvents the need to solve local optimization problems at each iteration, resulting in a \textit{single-loop} algorithm. Instead of applying the $\operatorname{argmin}$ operator, DGA uses a first-order approximation of each agent’s objective function based solely on local gradient information, thereby significantly reducing the per-iteration computational burden. To ensure feasibility, the proposed algorithm employs projection onto local constraint sets at each iteration, maintaining variable updates within the allowable operational limits. Another appealing feature of our algorithmis that it does not require exchanging sensitive information, such as the gradient and primal variables, making it privacy-preserving. Under mild and standard assumptions, we rigorously establish the convergence properties of the algorithm. Specifically, when the cost function is convex and Lipschitz-smooth, we establish a nonergodic convergence rate of $o(1/k)$. Under the additional assumption of strong convexity, the convergence further improves to a linear rate. Numerical simulation demonstrates that our algorithm outperforms state-of-the-art methods in terms of convergence speed, both with respect to the number of iterations and the overall computation time.

Our main contributions in this paper are summarized as follows.
\begin{enumerate}
    \item We propose a novel gradient-based distributed optimization algorithm for solving optimization problems with coupled affine equality constraints by adapting an approximate method of multipliers. In contrast to many existing approaches \cite{chang2014multi, chen2017admm, su2021distributed,falsone2023augmented, zhang2020distributed, wu2025distributed, zhu2025deed, yang2016distributed, shi2024distributed} that require solving a local optimization problem at each iteration, the proposed method eliminates the need for the use of $\operatorname{argmin}$ to solve local optimization problem per iteration. Instead, it relies solely on local gradient information and projection onto feasible sets, thereby significantly reducing per-iteration computational complexity. Unlike existing methods, the proposed algorithm does not require any special initialization and guarantees the feasibility of each agent’s local constraints throughout the iterations. 
    \item For convex cost functions, we prove that the KKT errors converge to $0$ at a nonergodic rate of $o\left(1/k\right)$. For strongly convex functions, we establish linear convergence. In contrast, existing algorithms based on similar first-order approximations, such as Mirror-PG-EXTRA \cite{nedic2018improved}, which lacks explicit convergence guarantees, and DuSPA \cite{xu2018dual}, which achieves a non-ergodic rate of $O(1/k)$ for convex case, offer weaker theoretical assurances. 
\end{enumerate}

The remainder of this paper is organized as follows. Section \ref{sec: Problem Formulation} introduces the problem formulation. Section \ref{sec: Algorithm Development} develops our single-loop distributed algorithm. Section \ref{sec: Convergence Analysis} establishes the convergence analysis. Section \ref{sec: Simulation} provides numerical experiments to validate our algorithm. Section \ref{sec: Conclusion} concludes the paper.

\textbf{Notation}: For a differentiable function $f: \mathbb{R}^{p} \rightarrow \mathbb{R}$, its gradient at $x \in \mathbb{R}^{p}$ is denoted by $\nabla f(x)$. 
For any, possibly non-differentiable, function $g: \mathbb{R}^{p} \rightarrow \mathbb{R}$, its subdifferential, i.e., its subgradient set at $x \in \mathbb{R}^{p}$, is represented by $\partial g(x)$. 
% A convex function \( f: X \rightarrow (-\infty, +\infty] \) is called proper if \( f(x) > -\infty \) for all \( x \in X \) and $f(x)$ is not trivially equal to $+\infty$. 
% The relative interior of a set $S$, $\operatorname{ri}(S)$, is defined as $\operatorname{ri}(S)=\left\{x \in S \mid \exists \epsilon>0\right.$ such that $\left.B_\epsilon(x) \cap \operatorname{aff}(S) \subseteq S\right\}$, where $\operatorname{aff}(S)$ is the affine hull of $S$ and $B_\epsilon(x)$ is a ball of radius $\epsilon$ centered on $x$.
The symbols $0_p$, $1_p$, and $I_p$ are used to denote the $p$-dimensional all-zero vector, all-one vector, and identity matrix, respectively. In addition, $O_{m\times p}$ denotes the $(m\times p)$-dimensional zero matrix. Let $\otimes$ denote the Kronecker product, $\langle\cdot, \cdot\rangle$ be the Euclidean inner product, and $\|\cdot\|$ be the $\ell_2$ norm. Given $A \in \mathbb{R}^{p \times p}, A^{\frac{1}{2}}$ is the square root of $A$ (i.e., $A^{\frac{1}{2}} A^{\frac{1}{2}}=A$ ), and we write $A \succeq O_{p \times p}$ if it is positive semidefinite and $A \succ O_{p \times p}$ if it is positive definite. For any $A \succeq O_{p \times p}$ and ${x} \in \mathbb{R}^p,\|{x}\|_A^2:={x}^{\rm{T}} A {x}$,
% the eigenvalues of $A$ are denoted by $\lambda_{p}(A)\geq \ldots \geq \lambda_{1}(A)$ in descending order and 
the spectral norm of $A$ is denoted by $\|A\|$. For any matrix $A$, $A^{\dagger}$ is $A^{\prime}$ 's pseudoinverse. In addition, we denote the block diagonal matrix by $\operatorname{diag}(A_1, A_2, \ldots, A_n)$ whose blocks are $A_1, A_2, \ldots, A_n$. The projection of a vector $x$ onto $X$ is denoted by $\mathcal{P}_{X}(x)$.

\section{Problem Formulation}\label{sec: Problem Formulation}
\subsection{Optimization Problem}
% Consider the EDP for a distributed power network consisting of multiple power generators. Local constraints exist for the output power of each individual generator. Besides, the total output power of all the generators needs to meet a certain demand of the power network. Suppose each generator has a local loss function $ f_i (x_i) $ and all the generators attempt to minimize the total loss while meeting total demands and satisfying individual limits, which are formulated as
% \begin{equation}\label{}
% \begin{aligned}
%     \min _{x} &~ f(x)= \sum_{i=1}^n f_i(x_i) \\
%     \operatorname{s.t.} &~  \sum_{i=1}^n A_i x_i = \sum_{i=1}^n d_i, \quad \underline{x}_i \leq x \leq \Bar{x}_i, \quad \forall i \in \mathcal{V},
% \end{aligned}
% \end{equation}
% where $x = (x_1,x_2,\ldots, x_n) \in \mathbb{R}$ represents the outpower of the generators, $\underline{x}_i$ and $\Bar{x}_i$ are lowest and highest limits of the generator $i \in \mathcal{V}$, respectively, and $d > 0$ is the power demand of the network.
This paper considers the following affine coupled-constraint problem,
\begin{align}\label{affine_constraint_problem}
    \min _{x \in X} f(x)= \sum_{i=1}^n f_i(x_i), \quad 
    \operatorname{s.t.}\sum_{i=1}^n A_i x_i = \sum_{i=1}^n d_i,
\end{align}
where $x=\left[ x_1^{\mathrm{T}}\, x_2^{\mathrm{T}}\, \ldots\, x_n^{\mathrm{T}}\right]^{\mathrm{T}} \in \mathbb{R}^{np}$, $x_i \in \mathbb{R}^p, i \in \mathcal{V}$, are the decision variables, $d_i \in \mathbb{R}^{m}$, $ A_i \in \mathbb{R}^{m \times p}$ and $\{A_i\}_{i=1}^n$ are not all zero, the local constraint for each agent is a convex constraint set $X_i \subseteq \mathbb{R}^p$. Accordingly, we define $X=X_1 \times X_2 \times \cdots \times X_n$. In problem \eqref{affine_constraint_problem}, we only require $m \leq np$, which is weaker than \cite{wu2025distributed}, where $A_i$ is required to have a full row rank.

% In modern distribution networks, energy management problems in power systems can be formulated in the form of problem \eqref{affine_constraint_problem} \cite{xie2024virtual}. Specifically, in the economic dispatch problem, local constraints exist for the output power of each individual generator. Besides, the total output power of all the generators needs to meet a certain demand of the power network. All the generators attempt to minimize the total loss while meeting total demands and satisfying individual limits. 
% On the other hand, in the context of communication systems, each transmitter or user terminal operates under its own physical limitations, such as maximum transmission power and bandwidth availability \cite{wang2024distributed}.
% At the same time, these individual decisions are inherently coupled through shared communication resources, interference among users, limited spectrum capacity, and network backhaul constraints, which together determine the system’s overall performance. The goal is to coordinate resource usage across the network so that communication efficiency is maximized and interference is mitigated, while ensuring that both local and global requirements are satisfied.

% In Sections \ref{subsec: Sublinear} and \ref{subsec: linear}, we analyze two cases: when $X_i \subset \mathbb{R}^p$ (i.e., constrained case), and when $X_i=\mathbb{R}^p$ (i.e., unconstrained case), for all $i \in \mathcal{V}$. 
\begin{assumption} \label{ass: smooth}
    For each $i \in \mathcal{V}$, the function $f_i$ is convex, differentiable and $l_{f}$-smooth, i.e.,
    $$
    \left\|\nabla f_i\left(x_i^{\prime}\right)-\nabla f_i\left(x_i\right)\right\| \leq l_{f}\left\|x_i^{\prime}-x_i\right\| \quad \forall x_i^{\prime}, x_i \in \mathbb{R}^p,
    $$
    for some $l_f > 0$.
\end{assumption}

\begin{assumption}[Slater's condition] \label{ass: slater's condition}
    The solution set is nonempty, and there exist some points in the relative interior of $X$, $\rm{ri}\{X\}$, such that the desired power constraint $\sum_{i=1}^n A_i x_{i} = \sum_{i=1}^n d_i$ is satisfied.   
\end{assumption}
\subsection{Communication Model}
We consider a networked system consisting of $n$ agents labeled by the index set $\mathcal{V}=\{1,2, \ldots, n\}$, where the communication structure among agents is modeled by a connected undirected graph $\mathcal{G}=(\mathcal{V}, \mathcal{E})$. Each vertex $i \in \mathcal{V}$ corresponds to an individual agent, and the edge set $\mathcal{E} \subseteq\{(i, j) \mid i \neq j, i, j \in \mathcal{V}\}$ encodes the bidirectional communication links between agents. An edge $(i, j) \in \mathcal{E}$ is associated with a positive scalar weight $p_{i j}=p_{j i}>0$, representing the interaction strength between agents $i$ and $j$. The communication topology is further characterized by the adjacency matrix $ P=\left[p_{i j}\right] \in \mathbb{R}^{n \times n}$, where $p_{i j}>$ 0 if and only if $(i, j) \in \mathcal{E}$, and $p_{i j}=0$ otherwise. For each agent $i \in \mathcal{V}$, we define its local neighborhood as $\mathcal{N}_i=\{j \in \mathcal{V} \mid(i, j) \in \mathcal{E}\}$. Communication is assumed to be local, i.e., each agent can exchange information only with its immediate neighbors.
% In addition, we assume that the undirected graph $G$ connected, and the weight matrix $P$ associated with the graph $G$ is doubly stochastic, i.e., $P 1_n = 1_n, ~1_n^{\rm{T}} P = 1_n^{\rm{T}}$.

% Proper convex function never takes on the value $-\infty$ and also is not identically equal to $\infty$.

\section{Algorithm Development} \label{sec: Algorithm Development}
By introducing the Lagrangian multiplier $y\in \mathbb{R}^p$, the Lagrangian function associated with problem \eqref{affine_constraint_problem} is defined as 
$f(x) + \left\langle \delta, \sum_{i=1}^n (A_i x_i - d_i) \right \rangle$.
The dual problem to problem \eqref{affine_constraint_problem} is then given as
\vspace{-0.3cm}
\begin{align} \label{dual}
    \underset{\delta \in \mathbb{R}^m}{\max} \underset{ x \in X }{\inf } L( x, \delta) \!
    % = & \ \underset{\delta \in \mathbb{R}^p}{\max} \ \underset{ x }{\inf } \left\{ F( x) + \left\langle \delta,  \sum_{i=1}^n (x_i - {\color{blue}d_i)} \right\rangle\right\} \notag\\
    % = & \ \underset{\delta \in \mathbb{R}^p}{\max} \ \underset{ x }{\inf } \left\{F( x) + \left\langle \delta,  \sum_{i=1}^n (x_i - {\color{blue}d_i)} \right\rangle \right\} \notag\\
    = &  \underset{\delta \in \mathbb{R}^m}{\max} \underset{ x \in X }{\inf } \! \left\{ \! \sum_{i=1}^n \left( f_i(x_i) \! + \! \langle \delta, A_i x_i \! - \! d_i \rangle  \right) \! \right\}\notag\\
    = &  \underset{\delta \in \mathbb{R}^m}{\max} \sum_{i=1}^n \underset{x_i \in X_i }{\inf }  \left\{ f_i(x_i)  + \! \langle \delta, A_i x_i \! -\! d_i  \rangle \right\}  \notag\\
    = &  \underset{\delta \in \mathbb{R}^m}{\max} \sum_{i=1}^n - f_i^*(-\delta) - \langle \delta, d_i \rangle ,
\end{align}  
where, for each $i \in \mathcal{V}$, the conjugate function $f_i^*(\delta)$ is defined as
% \begin{align*}
%     & f_i^*(\delta)\\
%     = &\ \underset{x_i \in X_i}{\sup} \langle \delta,x_i \rangle  - f_i(x_i) \\
%     = &\ -\underset{x_i \in X_i}{\inf} f_i(x_i) - \langle \delta,x_i \rangle .
% \end{align*}
\begin{align*}
     f_i^*(\delta)
    % = &\ \underset{x_i \in X_i}{\sup} \left\{\langle \delta,A_i x_i \rangle  - f_i(x_i) \right\} \\
    = -\underset{x_i \in X_i}{\inf} \left\{f_i(x_i) - \langle \delta,A_i x_i \rangle \right\}.
\end{align*}
Let $g_i(\delta) = f_i^*(-\delta) + \langle \delta, d_i \rangle$. Then, the dual problem in (\ref{dual}) can be written as %min _{y \in \mathbb{R}^p} \sum_{i=1}^n g_i(y).$
$\min _{\delta \in \mathbb{R}^p} \sum_{i=1}^n g_i(\delta)$. To solve this problem in a distributed manner, each agent $i \in \mathcal{V}$ maintains a local copy $y_i$ of the global variable $\delta$. The resulting formulation is a distributed consensus optimization problem
\begin{align}\label{consensusproblem}
    \min _{y \in \mathbb{R}^{n m} } G(y)=\sum_{i=1}^n g_i(y_i), \quad
    \operatorname{s.t.} ~ y_i = y_j, ~ i,j \in \mathcal{V}, 
\end{align}
where $y = \left[ y_1^{\mathrm{T}}\, y_2^{\mathrm{T}}\, \ldots\, y_n^{\mathrm{T}}\right]^{\mathrm{T}}$. In \cite{qiu2023stochastic}, it has been shown that, when the topology of the network is connected, problem \eqref{consensusproblem} is equivalent
to the following compact form:
% \begin{align}\label{compactdualoptimization}
%     \min _{y \in \mathbb{R}^{n m} }& ~ G(y)=\sum_{i=1}^n g_i(y_i), ~\\
%     \operatorname{s.t.}& ~ W^{\frac{1}{2}} y = 0_{nm}. \notag
% \end{align}
\begin{align}\label{compactdualoptimization}
    \min _{y \in \mathbb{R}^{n m} } G(y)=\sum_{i=1}^n g_i(y_i), \quad
    \operatorname{s.t.} W^{\frac{1}{2}} y = 0_{nm}. 
\end{align}
where $W = \mathcal{L} \otimes I_m \succeq {O}_{nm}$ and $\mathcal{L}$ is the Laplacian matrix of the network.
% $$
% [H]_{ij} = \left\{\begin{array}{ll}
% \sum_{s \in \mathcal{N}_i} H_{i s}, & i=j, \\
% -p_{i j}, & j \in \mathcal{N}_i, \\
% 0, & \text { otherwise, }
% \end{array} \quad \forall i, j \in \mathcal{V},\right.
% $$
By introducing the augmented Lagrangian function $L_{\rho}(y,v) = G(y)-v^{\rm{T}} W^{\frac{1}{2}} y+\frac{\rho}{2}\left\|W^{\frac{1}{2}} y\right\|^2, \rho>0$, we can linearize the method of multipliers \cite{boyd2011distributed} for solving \eqref{compactdualoptimization}, and the resulting update is given as follows,
\begin{subequations}
    \begin{alignat}{2}
        y^{k+1} =& \, \underset{y}{\operatorname{argmin}} \Big\{\langle \nabla G(y^k), y \rangle - \langle W^{\frac{1}{2}} v^k, y \rangle + \rho \langle W y^k, y \rangle \notag \\
        & \left.+ \frac{\eta}{2} \left\| y - y^k \right\|^2\right\} \notag \\
        =&\, y^k - \frac{1}{\eta}(\nabla G(y^k) + W^{\frac{1}{2}} v^k - \rho W y^k), \label{originialprimal}\\
        v^{k+1} =& \,v^k-\rho W^{\frac{1}{2}} y^{k+1},
    \end{alignat}
\end{subequations}
where $\nabla G(y^k) = \left[\nabla g_1^{\rm{T}}(y_1^k)\, \nabla g_2^{\rm{T}}(y_2^k)\, \ldots\, \nabla g_n^{\rm{T}}(y_n^k)\right]^{\rm{T}} \in \mathbb{R}^{nm}$ and 
\begin{align}\label{gradient g}
    \nabla g_i(y_i^k) 
    &=-\nabla f_i^*(-y_i^k)+d_i \notag \\
    &=-\underset{z \in X_i}{\operatorname{argmin}} ~\left\{f_i(z) + \langle A_i z, y_i^k\rangle \right\} +d_i \notag \\
    & = -A_i x_i^{k+1} + d_i.  
\end{align}
Let $A \in \mathbb{R}^{nm \times np}$ be a block-diagonal matrix defined as $A = \operatorname{diag}(A_1, A_2, \ldots, A_n)$ and $d =
[d_1^{\rm{T}}\,d_2^{\rm{T}}\,\ldots\,d_n^{\rm{T}}]^{\rm{T}}$. Then, $\nabla G(y^k)= A x^k + d$,
and the linearized method of multipliers for \eqref{compactdualoptimization} is given as 
\begin{subequations}
\begin{alignat}{3}
x^{k+1} = &\ \underset{z \in X}{\operatorname{argmin}} ~\left\{f(z) + \langle A z, y^k\rangle \right\},  \label{original x update}\\
y^{k+1} = &\ y^k - \frac{1}{\eta}(- A x^{k+1} + d- \lambda^k + \rho W y^k), \\
\lambda^{k+1}= &\ \lambda^k-\rho W y^{k+1} ,
\end{alignat}    
\end{subequations}
where $\lambda^k = W^{\frac{1}{2}} v^k$.
In the above algorithm, due to the presense of the local optimization problem at each iteration in \eqref{original x update} which might not be efficiently solved, we approximate the solution of $\operatorname{argmin}$ in \eqref{original x update} with the projection of $(x^k - \alpha (\nabla f(x^k) + A^{\rm{T}} y^k))$ onto set $X$. Substituting the projection into the dual ascent steps yields our proposed algorithm:
\begin{subequations}\label{all_update}
\begin{alignat}{3}
x^{k+1} = &\ \mathcal{P}_{X} (x^k - \alpha (\nabla f(x^k) + A^{\rm{T}} y^k)), \label{x update}\\
y^{k+1} = &\ y^k - \frac{1}{\eta}(-A x^{k+1} + d- \lambda^k + \rho W y^k), \label{y update}\\
\lambda^{k+1}= &\ \lambda^k-\rho W y^{k+1} \label{lambda update},
\end{alignat}    
\end{subequations}
where $\mathcal{P}_{X}$ is the projection operator over the convex set $X$.
% In both simulation and proof, it has been verified that the parameter $\eta$ and $\theta$ must satisfy $I - \frac{\theta}{\eta}W \succ 0$. 
The distributed implementation of \eqref{all_update} over the undirected network $(\mathcal{V}, \mathcal{E})$ is detailed in Algorithm \ref{algorithm}. 

\begin{algorithm}[ht]
\caption{Distributed Gradient-based Algorithm}
\label{algorithm}
\begin{algorithmic}[1]
\State \textbf{Initialization:} $x^0 = \left[ (x_1^0)^{\mathrm{T}}\, (x_2^0)^{\mathrm{T}}\, \ldots\, (x_n^0)^{\mathrm{T}}\right]^{\mathrm{T}}$, $y^0=\left[ (y_1^0)^{\mathrm{T}}\, (y_2^0)^{\mathrm{T}}\, \ldots\, (y_n^0)^{\mathrm{T}}\right]^{\mathrm{T}}$ are arbitrarily set. $\lambda^0=\left[ (\lambda_1^0)^{\mathrm{T}}\, (\lambda_2^0)^{\mathrm{T}}\, \ldots\, (\lambda_n^0)^{\mathrm{T}}\right]^{\mathrm{T}}$ is intialized to satisfy $(1^{\rm{T}}_{m} \otimes I_n) \lambda^0 = 0_m$ . Each agent $i \in \mathcal{V}$ sends the variable $y_i^{0}$ to its neighbors $j\in \mathcal{N}_i$. After receiving the information from its neighbors, each agent $i \in \mathcal{V}$ computes the aggregated information $t_i^{0} = \sum_{j \in \mathcal{N}_i} p_{i j}(y_i^{0}-y_j^{0})$. 
\For{$k = 0,1,2,\ldots, $} 
\State For each agent $i \in \mathcal{V}$ computes $x_i^{k+1} = \mathcal{P}_{X_i} (x_i^k - \alpha (\nabla f(x_i^k) + A_i^{\rm{T}} y_i^k))$.
\State Each agent $i \in \mathcal{V}$ updates $y_i^{k+1} = y_i^k - \frac{1}{\eta}(- A_i x_i^{k+1} + d_i- \lambda_i^k + \rho t_{i,k}) $ and sends the variable $y_i^k$ to its neighbors $j\in \mathcal{N}_i$.
\State Each agent $i \in \mathcal{V}$ sends the variable $y_i^{k+1}$ to its neighbors $j\in \mathcal{N}_i$. After receiving the information from its neighbors, each agent $i \in \mathcal{V}$ computes the aggregated information $t_i^{k+1} = \sum_{j \in \mathcal{N}_i} p_{i j}(y_i^{k+1}-y_j^{k+1})$. 
\State Each agent $i \in \mathcal{V}$ updates $\lambda_i^{k+1} = \lambda_i^k - t_i^{k+1}$.
\State Set $k \gets k+1$ and go to Step 3 until a certain stopping criterion is satisfied, e.g., maximum number of iterations.
\EndFor
\end{algorithmic}
\end{algorithm}

\section{Convergence Analysis}\label{sec: Convergence Analysis}
In this section, we consider two cases, with and without local constraints, i.e., $X_i \subset \mathbb{R}^p$ and $X_i = \mathbb{R}^p$, $i \in \mathcal{V}$. For these two cases, we establish sublinear and linear convergence rates, respectively.

We denote the optimal solution of problem \eqref{affine_constraint_problem} by $h^* = [(x^*)^{\mathrm{T}}\,(y^*)^{\mathrm{T}}\,(\lambda^*)^{\mathrm{T}}]^{\mathrm{T}}$, where  $x^*=$ $[(x_1^*)^{\mathrm{T}}\,(x_2^*)^{\mathrm{T}}\, \ldots \,\\(x_n^*)^{\mathrm{T}}]^{\mathrm{T}}$, $ y^*=[(y_1^*)^{\mathrm{T}}\,(y_2^*)^{\mathrm{T}}\, \ldots\,(y_n^*)^{\mathrm{T}}] ^{\mathrm{T}}$, and $\lambda^*=[(\lambda_1^*)^{\mathrm{T}}\,(\lambda_2^*)^{\mathrm{T}}\, \ldots\,(\lambda_n^*)^{\mathrm{T}}]^{\mathrm{T}}$. The optimal solution satisfies the first-order optimality conditions (KKT conditions) for convex constrained optimization problems as described in \cite{boyd2004convex}. In particular, the following conditions must hold. First, the dual variable $y^*$ lies in the null space of the $W$, i.e.,
\begin{equation}\label{consensus y^*}
    W y^* = 0_{nm},
\end{equation}
which indicates consensus of $y_i^*$ for all $i \in \mathcal{V}$. Second, let $\sigma_i(x_i)$ be the indicator function associated with the convex set $X_i$, $i \in \mathcal{V}$, i.e., $
\sigma_i\left(x_i\right)= 0 \text { if } x_i \in X_i$ and $ \sigma_i\left(x_i\right)= +\infty$, otherwise.
Then, the subdifferential condition associated with the objective and constraint functions is given by
\begin{equation}\label{partial sigma(x^*)}
    -A^{\rm{T}} y^* - \nabla f(x^*) \in \partial \sigma(x^*),
\end{equation}
where $\sigma(x) = \sum_{i=1}^n \sigma_i(x_i) $ and $\partial \sigma\left(x^*\right)$ denotes the subdifferential of the indicator function encoding local constraints at $x^*$. This reflects a standard KKT-type optimality condition under nonsmooth composite objectives. In addition, based on the definition of $L_{\rho} (y, v)$ and \eqref{gradient g}, every $v^*$ satisfying $W^{\frac{1}{2}} v^* = \nabla G(y^*)= -(A x^* - d)$ is an optimal dual variable for problem \eqref{compactdualoptimization}, and hence, 
\begin{equation}\label{opt_cond_lambda}
    \lambda^* = -(Ax^* - d).
\end{equation}
Since $\lambda^*$ is in the range space of $W$, \eqref{opt_cond_lambda} is equivalent to the coupled constraint, i.e., $(1^{\rm{T}}_{m} \otimes I_n) (A x^* - d) = 0_m.$
% \begin{equation}\label{opt_cond_feas}
%     (1^{\rm{T}}_{m} \otimes I_n) (A x^* - d) = 0_m.\end{equation}
% 
\subsection{Sublinear Convergence Rate}\label{subsec: Sublinear}
In this section, we carry out the convergence analysis of the proposed algorithm under Assumptions \ref{ass: smooth} and \ref{ass: slater's condition}. Denote the variable generated by the proposed algorithm at the $k$th iteration as $h^k = [(x^k)^{\rm{T}}\ (y^k)^{\rm{T}}\ (\lambda^k)^{\rm{T}}]^{\rm{T}}$. For $k \geq 0$, denote the difference of variables between the $(k+1)$th iteration and the $k$th iteration as $\Delta x^{k+1} = x^{k+1} - x^k$, $\Delta y^{k+1} = y^{k+1} - y^k$, $\Delta \lambda^{k+1} = \lambda^{k+1} - \lambda^k$, $\Delta \nabla f(x^{k+1})= \nabla f(x^{k+1}) - \nabla f(x^k)$ and $\Delta h^{k+1} = h^{k+1} - h^k $. The next two lemmas show the upper bound of the summation of $\| \Delta h^k \|^2_{\Omega}$ and the monotonicity of $\| \Delta h^k \|^2_{\Omega}$.
The relationship between two consecutive iterations are given in Lemma \ref{lem: h^k+1 -h}. 
\vspace{-0.2cm}
\begin{lem}\label{lem: h^k+1 -h}
   Let Assumptions \ref{ass: smooth} and \ref{ass: slater's condition} hold, If $\frac{\rho}{\eta}<\frac{1}{\lambda_{\max}(W)}$ and $\alpha< \frac{1}{l_f}$, we have the following inequalities which establish the monotonicity of $\| h^k-h^*\|_{{\Omega}}^2$ and provide an upper bound of the summation of the sequence $\Delta h^k$ generated by the proposed DGA, respectively, 
    \begin{align}\label{lemma1}
        &\, \| h^k-h^*\|_{{\Omega}}^2-\| h^{k+1} -h^*\|_{{\Omega}}^2 \notag \\
        \geq &\, \| h^{k+1} - h^k \|_{{\Omega}}^2 - \alpha l_f \| x^{k+1} - x^k \|^2 , ~k \geq 0,
    \end{align}
    \begin{equation}\label{eqlem: upperbounded}
        \sum_{t=0}^k\left\|\Delta h^{t+1}\right\|_{\Omega}^2<\frac{\| h^0-h^*\|_{{\Omega}}^2}{(1-\alpha l_f)}, ~ k \geq 1,
    \end{equation}
    where $\Omega \in \mathbb{R}^{(np+2nm)\times (np+2nm)}$ is defined as $$\Omega = \left[\begin{array}{ccc}
    I_{np} & -\alpha A^{\rm{T}} & O_{np\times nm} \\
    O_{nm\times np} & \alpha \eta (I_{nm} \!-\! \frac{\rho}{\eta} W) & O_{nm\times nm} \\
    O_{nm\times np} & O_{nm\times nm} & \frac{\alpha}{\rho}W^{\dagger}
    \end{array}\right].$$
\end{lem}
\vspace{-0.4cm}
\begin{pf}
By \eqref{y update} and \eqref{consensus y^*}, we have
\begin{equation}\label{lambdak-lambda*}
    \lambda^{k+1}-\lambda^*=\lambda^k-\lambda^*-\rho W\left(y^{k+1}-y^*\right).
\end{equation}
In addition, \eqref{lambda update} and \eqref{opt_cond_lambda} yield
\begin{align}\label{yk+1-y*}
    y^{k+1} - y^* = & \left(I_{nm} \!-\! \frac{\rho}{\eta} W\right) (y^k - y^*) \displaybreak[0] \notag \\
    & - \frac{1}{\eta}\left(A (x^{k+1} - x^*) - (\lambda^k - \lambda^*)\right).
\end{align}
Note that $(I_{nm} \!-\! \tfrac{\rho}{\eta} W) \succ O_{nm}$. Combining the above two equations leads to 
\begin{align}\label{relation_x_y_lambda}
    &\, A(x^{k+1} - x^*)\notag \\
    =&\, \eta \left(I_{nm} \!-\! \frac{\rho}{\eta} W\right) (y^{k+1} - y^k) -  (\lambda^{k+1} - \lambda^*).
\end{align}
In addition, by the optimal condition of \eqref{x update}, we have
\begin{equation}\label{opt_condition_x_update}
    -x^{k+1}+(x^k-\alpha( \nabla f(x^k)+A^{\rm{T}} y^k)) \in \partial \sigma\left(x^{k+1}\right).
\end{equation}
The convexity of the indicator function $\sigma$ leads to $\langle \partial \sigma\left(x^{k+1}\right) - \partial \sigma\left(x^*\right), x^{k+1}-x^* \rangle \geq 0.$ Therefore, in view of \eqref{partial sigma(x^*)} and \eqref{opt_condition_x_update}, we have
\begin{align}\label{convexity_sigma}
    0 \leq &\, \langle x^k-x^{k+1} -\alpha( (\nabla f\left(x^k\right)+A^{\rm{T}} y^k) \notag\\
    &  - ( A^{\rm{T}} y^*+\nabla f(x^*)) ), x^{k+1}-x^*\rangle.
\end{align}
Then, we can embed \eqref{relation_x_y_lambda} into \eqref{convexity_sigma} to obtain
\begin{align}\label{monotonicity 2}
    & -\alpha\langle \nabla f(x^k) - \nabla f(x^*), x^{k+1} - x^*\rangle \notag \displaybreak[0] \\
    \geq& \, \langle x^{k+1} - x^k, x^{k+1} - x^*\rangle \! + \! \alpha\langle A^{\rm{T}}( y^k - y^*), x^{k+1} - x^*\rangle \notag \displaybreak[0] \\
    = &\, \langle x^{k+1} - x^k, x^{k+1} - x^*\rangle \! + \! \alpha\langle A^{\rm{T}}(y^k \! - \! y^{k+1}), x^{k+1} - x^*\rangle \notag \\
    & + \alpha\langle A^{\rm{T}}(y^{k+1} - y^*), x^{k+1} - x^*\rangle \notag \displaybreak[0] \\
    = &\, \langle x^{k+1} - x^k, x^{k+1} - x^*\rangle \! + \! \alpha\langle A^{\rm{T}}(y^k \! - \! y^{k+1}), x^{k+1} - x^*\rangle \notag \\
    & + \alpha\left\langle y^{k+1} - y^*, \eta \left(I_{nm} \!-\! \frac{\rho}{\eta} W\right) (y^{k+1} - y^k) \right \rangle \notag \\
    & - \alpha\langle y^{k+1} - y^*,\lambda^{k+1} - \lambda^* \rangle .
\end{align}
By \eqref{lambda update} and \eqref{consensus y^*}, we have $\lambda^{k+1} - \lambda^k = -\rho W (y^{k+1} - y^*).$
% \begin{equation*}
%     \lambda^{k+1} - \lambda^k = -\rho W (y^{k+1} - y^*).
% \end{equation*}
Since $W$ is a semi-definite matrix with only one zero eigenvalue, we have the pseudo inverse matrix of $W$ is also semi-definite with only one zero eigenvalue. Thus, we have 
\begin{equation}
     -\frac{1}{\rho} W^{\dagger} (\lambda^{k+1} - \lambda^k) = y^{k+1} - y^*,
\end{equation}
and hence, $\langle y^* \!-\! y^{k+1} , \lambda^{k+1} \!-\! \lambda^* \rangle 
     \!=\!  \langle \frac{1}{\rho} W^{\dagger} (\lambda^{k+1} \!-\! \lambda^k),\lambda^{k+1} \!-\! \lambda^* \rangle,$
% \begin{align*}
%      \langle y^* \!-\! y^{k+1} , \lambda^{k+1} \!-\! \lambda^* \rangle 
%      \!=\! \left \langle \frac{1}{\rho} W^{\dagger} (\lambda^{k+1} \!-\! \lambda^k),\lambda^{k+1} \!-\! \lambda^* \right \rangle,
% \end{align*}
substitution of which into \eqref{monotonicity 2} yields
\begin{align}\label{monotonicity 3}
    & -\alpha\langle \nabla f(x^k) - \nabla f(x^*), x^{k+1} - x^*\rangle \notag \displaybreak[0] \\
    \geq &\, \langle x^{k+1} - x^k, x^{k+1} - x^*\rangle \! + \! \alpha\langle A^{\rm{T}}(y^k \! - \! y^{k+1}), x^{k+1} - x^*\rangle \notag \displaybreak[0] \\
    & + \alpha \left\langle y^{k+1} - y^*, \eta \left(I_{nm} \!-\! \frac{\rho}{\eta} W\right) (y^{k+1} - y^k) \right\rangle \notag \displaybreak[0] \\
    & + \frac{\alpha}{\rho} \langle W^{\dagger} (\lambda^{k+1} - \lambda^k),\lambda^{k+1} - \lambda^* \rangle \notag \displaybreak[0] \\
    = &\, \langle h^{k+1} - h^k, h^{k+1} - h^*\rangle_{\Omega} \notag \displaybreak[0] \\
    = & \frac{1}{2}(\| h^{k+1}\! -\!h^*\!\left\|_{{\Omega}}^2 - \right\| h^k\! -\!h^*\!\left\|_{{\Omega}}^2+\right\| \! h^{k+1} \! -\! h^k \|_{{\Omega}}^2).
\end{align}
By the Lipschitz smoothness and convexity of $f$, we have
% \begin{equation}\label{fxk+1 fxk + norm}
%     f(x^{k+1}) \leq f(x^k) + \langle \nabla f(x^k), x^{k+1} - x^k \rangle + \frac{l_f}{2} \| x^{k+1} - x^k \|^2,
% \end{equation}
% \begin{equation}\label{convex_fxk+1}
%     f(x^{k+1}) \geq f(x^*) + \langle \nabla f(x^*), x^{k+1} - x^* \rangle.
% \end{equation}
\begin{equation*}
    f(x^{k+1}) \leq f(x^k) + \langle \nabla f(x^k), x^{k+1} - x^k \rangle + \frac{l_f}{2} \| x^{k+1} - x^k \|^2,
\end{equation*}
\begin{equation*}
    f(x^{k+1}) \geq f(x^*) + \langle \nabla f(x^*), x^{k+1} - x^* \rangle.
\end{equation*}
Combining the above two inequalities, we obtain
\begin{align*}
    & f(x^k) - f(x^*) + \frac{l_f}{2} \| x^{k+1} - x^k \|^2 \notag \\
    \geq & \langle \nabla f(x^*), x^{k+1} - x^* \rangle - \langle \nabla f(x^k), x^{k+1} - x^k \rangle,
\end{align*}
which together with $f(x^k) - f(x^*) \leq \langle \nabla f(x^k), x^k - x^* \rangle$, results in 
\begin{align*}
    & \left\langle\nabla f\left({x}_k\right)-\nabla f\left({x}^{*}\right), {x}_{k+1}-{x}^{*}\right\rangle \geq  -\frac{l_f}{2} \| x^{k+1} - x^k \|^2,
\end{align*}
substitution of which into \eqref{monotonicity 3} results in \eqref{lemma1}.
% \begin{align*}
% &\, \alpha l_f \| x^{k+1} - x^k \|^2   \notag \\
% \geq &\, \| h^{k+1} -h^*\|_{{\Omega}}^2 - \| h^k-h^*\|_{{\Omega}}^2 + \| h^{k+1} - h^k \|_{{\Omega}}^2.
% \end{align*}

Noting that $\alpha< \frac{1}{l_f}$, we have
\begin{align}\label{v^k-v^{k+1} 2}
    \left\|\Delta h^{k+1}\right\|_{\Omega-\alpha l_f {\Omega}^{\prime}}^2
    \geq \left(1-\alpha l_f\right)\left\|\Delta h^{k+1}\right\|_{\Omega}^2 \geq 0,
\end{align}
where $\Omega^{\prime} = {\rm{diag}} \{ I_{np}, O_{nm}, O_{nm} \}$ and the second equality holds at $\Delta h^{k+1} = 0$.
Based on \eqref{v^k-v^{k+1} 2}, summing \eqref{lemma1} from $0$ to $k$ leads to \eqref{eqlem: upperbounded}.  \hfill $\Box$
\end{pf}
\begin{lem}\label{lem: nonincreasing}
    Let Assumptions \ref{ass: smooth} and \ref{ass: slater's condition} hold. If $\frac{\rho}{\eta}<\frac{1}{\lambda_{\max}(W)}$ and $\alpha< \frac{1}{l_f}$, then the sequence $\| \Delta h^k \|_{\Omega}$ generated by the proposed DGA algorithm is monotonically nonincreasing, that is,
    \begin{equation}
        \| \Delta h^k \|^2_{\Omega} \geq \| \Delta h^{k+1} \|^2_{\Omega}.
    \end{equation}
\end{lem}
\begin{pf}
    Similar to \eqref{relation_x_y_lambda}, we have the relationship among $\Delta y^{k+1},~\Delta x^{k+1}$ and $ \Delta \lambda^{k+1}$ as follows
    \begin{equation}\label{relation_deltayxlambda}
        \eta \left(I_{nm} \!-\! \frac{\rho}{\eta} W\right)(\Delta y^{k+1} - \Delta y^k) = A \Delta x^{k+1} -\Delta \lambda^{k+1}.
    \end{equation}
    Taking the difference of the $k$th update and the $(k+1)$th update of \eqref{opt_condition_x_update} yields 
    $$
    - (\Delta x^{k+1}- \Delta x^k)-\alpha ( \Delta\nabla f\left(x^k\right)+A^{\rm{T}} \Delta y^k)) \in \sigma_{k+1} - \sigma_k,
    $$
    which, in view of the convexity of $\sigma$, i.e., $\langle \Delta x^{k+1}, \sigma_{k+1} - \sigma_k \rangle \geq 0$, results in
    \begin{align}\label{resultof_optcondition}
        & \langle \Delta x^{k+1},  - \Delta x^{k+1} + \Delta x^k \rangle + \alpha \langle \Delta x^{k+1}, -  \Delta\nabla f\left(x^k\right) \rangle \notag \\
        & + \alpha \langle \Delta x^{k+1}, - A^{\rm{T}} \Delta y^k \rangle \geq 0.
    \end{align}
    The equivalent expression of the last term on the right hand side of the above equation is
    \vspace{-0.1cm}
    \begin{align}\label{last_term_above}
        &\ \alpha \langle \Delta x^{k+1}, - A^{\rm{T}} \Delta y^k \rangle \displaybreak[0] \notag \\ 
        = &\ \alpha \langle \Delta x^{k+1}, A^{\rm{T}} (\Delta y^{k+1} - \Delta y^k)\rangle - \alpha \langle \Delta x^{k+1}, - A^{\rm{T}} \Delta y^{k+1} \rangle \displaybreak[0] \notag \\ 
        = &\ \alpha \langle \Delta x^{k+1}, A^{\rm{T}} (\Delta y^{k+1} - \Delta y^k)\rangle \displaybreak[0] \notag \\
        & - \alpha \left \langle \eta \left(I_{nm} \!-\! \frac{\rho}{\eta} W\right)(\Delta y^{k+1} - \Delta y^k), - \Delta y^{k+1} \right \rangle  \displaybreak[0] \notag \\ 
        & - \alpha \langle \Delta \lambda^{k+1}, - \Delta y^{k+1} \rangle \displaybreak[0] \notag \\ 
        = &\ \alpha \langle \Delta x^{k+1}, A^{\rm{T}} (\Delta y^{k+1} - \Delta y^k)\rangle \displaybreak[0] \notag \\
        & - \alpha \left\langle \eta \left(I_{nm} \!-\! \frac{\rho}{\eta} W\right)(\Delta y^{k+1} - \Delta y^k), - \Delta y^{k+1} \right \rangle \displaybreak[0] \notag \\ 
        & - \alpha \left\langle \Delta \lambda^{k+1}, \frac{1}{\rho} W^{\dagger} (\Delta \lambda^{k+1} - \Delta \lambda^k) \right\rangle,
    \end{align}
    where the second equality results from \eqref{relation_deltayxlambda} and the last equality is due to \eqref{lambda update}.
    Substituting \eqref{last_term_above} back into \eqref{resultof_optcondition}, we have
    \vspace{-0.1cm}
    \begin{align}\label{nonincreasing 2}
    & -\alpha\langle \Delta \nabla f(x^k), \Delta x^{k+1} \rangle \notag \\
    \geq&\, \langle \Delta x^{k+1} \! - \! \Delta x^k, \Delta x^{k+1} \rangle \! + \! \alpha\langle A^{\rm{T}}\!(\!\Delta y^k \! - \! \Delta y^{k+1}\!), \Delta x^{k+1}\rangle \notag \\
    & + \alpha \left \langle \eta \left(I_{nm} \!-\! \frac{\rho}{\eta} W\right) (\Delta y^{k+1} - \Delta y^k), \Delta y^{k+1} \right\rangle \notag \\
    & + \alpha\langle \Delta \lambda^{k+1} -\Delta \lambda^k,\Delta \lambda^{k+1} \rangle \notag \\
    =&\, \frac{1}{2}(\| \Delta h^{k+1} \|_{\Omega}^2 - \| \Delta h^k \|_{\Omega}^2 + \|\Delta h^{k+1} - \Delta h^k \|_{\Omega}^2).
    \end{align}
    On the other hand
    \begin{align}\label{scale_Delta_nabla_f}
        & -\alpha\langle \Delta \nabla f(x^k), \Delta x^{k+1} \rangle \notag \\
        = & -\alpha\langle \Delta \nabla f(x^k), \Delta x^{k+1} - \Delta x^k \rangle -\alpha\langle \Delta \nabla f(x^k), \Delta x^k \rangle \notag \\
        \leq &\, \left(\frac{\alpha}{2 \xi}- \frac{\alpha}{l_f}\right) \| \Delta \nabla f(x^k) \|^2 + \frac{\alpha \xi}{2}\| \Delta x^{k+1} - \Delta x^k \|^2,  
    \end{align}
    where for any $\xi > 0$ the last inequality is due to the triangular inequality and the Lipschitz smoothness of $f$. Taking $\xi > l_f/2$ and combining \eqref{scale_Delta_nabla_f} and \eqref{nonincreasing 2} results in
    \begin{align*}
        & \| \Delta h^k \|_{\Omega}^2 \notag \\
        \geq & \| \Delta h^{k+1} \|_{\Omega}^2 + \|\Delta h^{k+1} - \Delta h^k \|_{\Omega}^2 - \alpha \xi\| \Delta x^{k+1} - \Delta x^k \|^2,
    \end{align*}
    which, in view of $\xi \in (\frac{l_f}{2}, \frac{1}{\alpha})$, completes the proof of Lemma \ref{lem: nonincreasing}. \hfill $\Box$
\end{pf}

\begin{thm}\label{the: subl inear convergence}
     Let Assumptions \ref{ass: smooth} and \ref{ass: slater's condition} hold. If $\frac{\rho}{\eta}<\frac{1}{\lambda_{\max}(W)}$ and $\alpha \in \left(0,~ \min \left\{ \frac{1}{l_f}, \frac{4\lambda_{\min}(\eta I_{nm} - \rho W)}{\|A\|_2^2} \right \} \right)$, then $\|\Delta h^k\|_{\Omega}^2 = o\left(1/k\right)$, and, consequently, the errors in the KKT conditions vanish at the rate of $o\left(1/k\right)$,
     \begin{subequations}
     \begin{alignat}{2}
         &\, \left \| \sum_{i=1}^n A_i x_i^{k+1} - d_i \right \|^2 = o\left(\frac{1}{k}\right), \label{feasibilityerror} \\
         &\,\|\alpha( \nabla f(x^{k+1})+A^{\rm{T}} y^k)) + \nu^{k+1}\|^2 = o\left(\frac{1}{k}\right), \label{firstorderopterror}
     \end{alignat}
     \end{subequations}
     where $\nu^{k+1} \in \partial \sigma\left(x^{k+1}\right)$.
     % $h^k$ converges to an optimal solution $h^*$ of problem \eqref{affine_constraint_problem}. 
\end{thm}
\vspace{-0.2cm}
\begin{pf}
    Combining Lemmas \ref{lem: h^k+1 -h} and \ref{lem: nonincreasing} with \cite[Proposition~1]{nedic2018improved}, we have that the sequence $\Delta h^k$ satisfies the nonergodic rate, i.e., $\|\Delta h^k\|_{\Omega}^2 = o(1/k)$.
    Let $S = \frac{1}{2}(\Omega + \Omega^{\rm{T}})$ and $\|\Delta h^k\|^2_{\Omega} = \|\Delta h^k\|^2_{S}$. Since 
    $ 
    \alpha < \frac{4\lambda_{\min}(\eta I_{nm} - \rho W)}{\|A\|_2^2},
    $
    we have $\alpha\eta\!\left(I_{nm} - \tfrac{\rho}{\eta}W\right)
    - \tfrac{\alpha^{2}}{4}A A^{\!\top} \succeq 0$ and $S$ is thus positive semidefinite, 
    which directly implies that
    \begin{subequations}
    \begin{alignat}{3}
        &\, \| x^{k+1} - x^k \|^2 = o\left(\frac{1}{k}\right), \label{x_o1k} \\
        &\, \| y^{k+1} - y^k \|^2_{\eta \left(I_{nm} \!-\! \frac{\rho}{\eta} W\right)} = o\left(\frac{1}{k}\right). \label{y_o1k}
        % &\, \| \lambda^{k+1} - \lambda^k \|^2_{W^{\dagger}}=o\left(\frac{1}{k}\right). \label{lambda_o1k}
    \end{alignat}
    \end{subequations}
    % Combining \eqref{lambda_o1k} and \eqref{lambda update} yields the convergence rate of the consensus error of $y^{k+1}$, i.e.,
    % \begin{equation}
    %     \| y^{k+1} \|_W^2 = o\left(\frac{1}{k}\right).
    % \end{equation}
    Since both $\lambda^k$ and $\lambda^*$ lie in the null space of $W$, left-multiplying both sides of \eqref{relation_x_y_lambda} by $(1^{\rm{T}}_{m} \otimes I_n)$ gives
    \begin{align}
        &\, \eta (1^{\rm{T}}_{m} \otimes I_n) \left(I_{nm} \!-\! \frac{\rho}{\eta} W\right) (y^{k+1} - y^k)\notag \\
        =&\, (1^{\rm{T}}_{m} \otimes I_n)(A x^{k+1} - d) = \sum_{i=1}^n A_i x_i - d_i \notag .
    \end{align}
    Taking norms of both sides and using~\eqref{y_o1k}, we obtain 
    \begin{align}
        & \left \| \sum_{i=1}^n A_i x_i^{k+1} - d_i \right \|^2 \notag \\
        =&\, \eta \left\|1^{\rm{T}}_{m} \otimes I_n\right\|^2 \| (I_{nm} \!-\! (\rho/\eta) W)\| \| y^{k+1} - y^k \|^2_{\eta \left(I_{nm} \!-\! \frac{\rho}{\eta} W\right)}  \notag ,
        % \\=&\, o\left(\frac{1}{k}\right).
    \end{align}
    which establishes~\eqref{feasibilityerror}.
    By \eqref{opt_condition_x_update}, there exists $\nu^{k+1} \in \partial \sigma\left(x^{k+1}\right)$ such that
    $$
    \| x^{k+1} - x^k \| = \|\alpha( \nabla f(x^k)+A^{\rm{T}} y^k)) + \nu^{k+1}\|,
    $$
   Using the Lipschitz continuity of $\nabla f$ and inequality~\eqref{x_o1k} yields
    \begin{align}
        &\, \|\alpha( \nabla \! f(x^{k+1})\!+\!A^{\rm{T}} y^k)) \!+\! \nu^{k+1}\| \notag \\
        \leq &\, \|\alpha( \nabla \! f(x^k)\!+\!A^{\rm{T}} y^k)) \!+\! \nu^{k+1}\|   \!+\! \|\nabla \! f(x^{k+1}) \!-\! \nabla \! f(x^k)\| \notag \\
        \leq &\, (\alpha l_f + 1) \| x^{k+1} - x^k \|, \notag 
    \end{align}
    which establishes~\eqref{firstorderopterror}.    \hfill $\Box$
    % From Lemma \ref{lem: h^k+1 -h}, it follows that the sequence $\{\| h^k - h^* \|^2_{\Omega}\}$ is bounded, which implies that $\{h^k\}$ is also bounded. 
    % Hence, the sequence $\{h^k\}_{k\ge0}$ admits at least one cluster point, denoted by $h_\infty$, 
    % and there exists a corresponding convergent subsequence $\{h_{k_i}\}_{i\ge0}$ such that $h_{k_i}\to h_\infty$. 

    % Taking the limit of the iterative relation \eqref{all_update} along this subsequence yields that $h_\infty$ satisfies the optimality conditions of the original problem, and thus $h_\infty$ is an optimal solution. 
    % Since $\{\| h^k - h^* \|^2_{\Omega}\}$ converges and is contractive with respect to the optimal solution set $S$, it follows that $h_\infty$ is the unique cluster point of $\{h^k\}$ . 
    % Consequently, the entire sequence $\{h^k\}_{k\ge0}$ converges to a unique point $h^*\in S$. 
    % Therefore, we conclude that $\{h^k\}$ converges to the optimal solution $h^*$ 
    % with a non-ergodic convergence rate of $o(1/k)$.
\end{pf}
% Theorem \ref{the: sublinear convergence} establishes a nonergodic sublinear convergence rate for the proposed distributed algorithm under standard smoothness and constraint qualification assumptions. 
% In particular, it shows that the sequence of iterates satisfies $\|\Delta h^k\|_{\Omega}^2 = o(1/k)$, 
% implying that the successive updates diminish asymptotically at a rate faster than $1/k$. 
% More importantly, this result directly guarantees the vanishing of the KKT residuals---both the feasibility error and the first-order optimality error---at the same asymptotic rate. 
% Specifically, the primal feasibility violation $\|\sum_i A_i x_i^{k+1} - d_i\|$ 
% and the stationarity residual $\|\alpha(\nabla f(x^{k+1}) + A^{\top} y^k) + \nu^{k+1}\|$ 
% converge to zero with an $o(1/k)$ decay rate. 
% This ensures that the algorithm progressively enforces both feasibility and first-order optimality conditions across the network. 
% The nonergodic rate result is particularly notable as it characterizes the \emph{pointwise} convergence behavior---as opposed to the averaged (ergodic) convergence rate---which represents a stronger property in distributed convex optimization. 
% Consequently, Theorem~\ref{the: sublinear convergence} provides a rigorous theoretical foundation for the asymptotic accuracy and stability of the proposed algorithm in solving linearly constrained convex optimization problems.
\vspace{-0.2cm}
\subsection{Linear Convergence Rate}\label{subsec: linear}
In this section, we consider the case of $X = \mathbb{R}^{np}$ and establish a linear convergence rate of Algorithm \ref{algorithm} under certain conditions. 
\vspace{-0.2cm}
\begin{assumption}\label{ass: strong convexity}
    For each $i \in \mathcal{V}$, the cost function $f_i$ is $\mu$-strongly convex, i.e.,
    $$
    \left\|\nabla f_i\left(x_i^{\prime}\right)-\nabla f_i\left(x_i\right)\right\| \geq \mu\left\|x_i^{\prime}-x_i\right\| \quad \forall x_i^{\prime}, x_i \in \mathbb{R}^p,
    $$
    for some $\mu > 0$.
\end{assumption}
\vspace{-0.2cm}
Since the constraint in problem \eqref{affine_constraint_problem} is affine and $X$ is convex, the feasible set is convex. Under a $\mu$-strongly convex objective function and Assumption \ref{ass: slater's condition}, the problem has a unique optimal solution in the feasible convex set.
% \begin{lemma}\label{lem: h^{k+1} - p^*<h^k - p^*-}
% Suppose that $X = \mathbb{R}^{np}$ and Assumptions \ref{ass: smooth}-\ref{ass: strong convexity} hold, for any $k>0$, we have
% \begin{equation}
% \begin{aligned}
% & \| h^{k+1} - p^* \|_{\Theta_1}^2 \\
% \leq & \| h^k - p^* \|_{\Theta_2}^2 - \left\|\Delta h^{k+1}\right\|_{\Omega}^2+\alpha \theta_0\left\|\Delta x^{k+1}\right\|^2 \\
% & -\left(\frac{2 \alpha}{l_f+\mu}-\frac{\alpha}{\theta_0} - \frac{\alpha^2}{\eta \theta_2}\right)\left\|\nabla f\left(x^k\right)-\nabla f\left(x^{*}\right)\right\|^2,
% \end{aligned}
% \end{equation}
% where $$\Theta_1 = \left[\begin{array}{ccc}
% \rho W - \frac{\rho^2}{\eta} W^2 & 0 & 0 \\
% 0 & (\alpha \eta I + \rho W)(I_{nm} \!-\! \frac{\rho}{\eta} W) & 0 \\
% 0 & 0 & \frac{\alpha}{\rho}W^{\dagger} + \frac{1}{\eta}I
% \end{array}\right]$$, $$\Theta_2 = \left[\begin{array}{ccc}
% (1+\frac{1}{\eta \theta_1} - \frac{2 \mu l_f \alpha}{l_f+\mu})I & 0 & 0 \\
% 0 & \alpha \eta (I_{nm} \!-\! \frac{\rho}{\eta} W) + A & 0 \\
% 0 & 0 & \frac{\alpha}{\rho}W^{\dagger} + \frac{1}{\eta}I - \frac{\rho}{\eta^2}W
% \end{array}\right].$$
% \end{lemma}
\vspace{-0.1cm}
\begin{thm}\label{the: linear convergence}
Suppose that $X = \mathbb{R}^{np}$ and Assumptions \ref{ass: smooth}-\ref{ass: strong convexity} hold. Suppose that the parameters satisfy $$
\eta \! > \! \max \! \left\{\! \frac{\alpha \|A\|^2 }{\left(\frac{2}{l_f \!+\! \mu}\!-\!\frac{1}{\theta_0}\right) \theta_2}, \frac{l_f \!+\! \mu}{2 \mu l_f \alpha} \!\left(\!\rho \| A \|^2 \| W \|\! +\! \frac{ \| A \|^2}{\theta_1} \!\right) \! \right\},
$$ $\alpha < \frac{2}{l_f + \mu}$, $\theta_0 \in (\frac{l_f + \mu}{2}, \frac{1}{\alpha}]$, and $\theta_1 + \theta_2 < 1$, the sequence $\left\{h^k\right\}_{k \geq 0}$ generated by the proposed algorithm will converge linearly to the optimal solution $h^*$ of the problem \eqref{affine_constraint_problem}. In particular, we have
\begin{equation}\label{theorem 2}
    \left\|h^{k+1}-h^*\right\|_{\Theta}^2 \leq \delta \left\|h^k-h^*\right\|_{\Theta}^2,
\end{equation}
where for some positive definite $\Theta$ and $\delta \in (0,1).$
\end{thm}
\vspace{-0.6cm}
\begin{pf}
We first establish the basic inequality showing the decreasing property of the optimality gap $\left\|h^k-h^*\right\|_{\Omega}$ based on the smoothness and strong convexity of $f$. By \eqref{monotonicity 3}, we have
\begin{align}\label{h-h*}
& \left\|h^{k+1}-h^*\right\|_{\Omega}^2-\left\|h^k-h^*\right\|_{\Omega}^2 \notag \displaybreak[0]\\
\leq & -\left\|\Delta h^{k+1}\right\|_{\Omega}^2-2 \alpha\left\langle\nabla f\left(x^k\right)-\nabla f\left(x^{*}\right), x^{k+1}-x^{*}\right\rangle \notag \displaybreak[0] \\
\leq& -\left\|\Delta h^{k+1}\right\|_{\Omega}^2-2 \alpha\left\langle\nabla f\left(x^k\right)-\nabla f\left(x^{*}\right), x^k-x^{*}\right\rangle \notag \displaybreak[0] \\
& -2 \alpha\left\langle\nabla f\left(x^k\right)-\nabla f\left(x^{*}\right), \Delta x^{k+1}\right\rangle .
\end{align}    
For the last term on the right-hand side of the second inequality in \eqref{h-h*}, using Young's inequality, for any given $\theta_0 > 0$ we have
\begin{align}\label{nabla f geq}
& \left\langle\nabla f\left(x^k\right)-\nabla f\left(x^{*}\right), \Delta x^{k+1}\right\rangle \notag \\
\geq &-\frac{\left\|\nabla f\left(x^k\right)-\nabla f\left(x^{*}\right)\right\|^2}{2 \theta_0}-\frac{\theta_0\left\|\Delta x^{k+1}\right\|^2}{2} .
\end{align}
Since $f$ is $l_f$-smooth and $\mu$-strongly convex, for the second term of \eqref{h-h*} we have
\begin{align}\label{nable f geq 2}
& \! \left\langle\nabla f\left(x^k\right) \! - \!  \nabla f\left(x^{*}\right), x^k\! - \! x^{*}\right\rangle \notag \\
\geq & \frac{1}{l_f\! + \! \mu}\!\left\|  \nabla f\left(x^k\right)\! - \! \nabla f\left(x^{*}\right)  \right\|^2\! \!+ \! \frac{\mu l_f}{l_f\! + \! \mu}\left\|x^k\! - \! x^{*}\right\|^2 .
\end{align}
Using \eqref{nabla f geq} and \eqref{nable f geq 2} in \eqref{h-h*}, we have
\begin{align}\label{omega_norm_1}
& \left\|h^{k+1}-h^*\right\|_{\Omega}^2-\left\|h^k-h^*\right\|_{\Omega}^2 \notag \\
\leq & -\left\|\Delta h^{k+1}\right\|_{\Omega}^2+\alpha \theta_0\left\|\Delta x^{k+1}\right\|^2-\frac{2 \mu l_f \alpha}{l_f+\mu}\left\|x^k-x^{*}\right\|^2 \notag  \\
& -\left(\frac{2 \alpha}{l_f+\mu}-\frac{\alpha}{\theta_0}\right)\left\|\nabla f\left(x^k\right)-\nabla f\left(x^{*}\right)\right\|^2.
\end{align}
Taking square of both sides of \eqref{lambdak-lambda*} yields
\begin{align}\label{based_on_lambdaupdate}
    \| \lambda^{k+1} - \lambda^*\|^2 = &\, \| \lambda^k - \lambda^* -\rho W (y^{k+1} - y^*) \|^2 \notag \\
    = &\,  \| \lambda^k - \lambda^* \|^2 + \rho^2 \| y^{k+1} - y^* \|_{W^2}^2 \notag  \\
    & - 2 \rho\langle W(y^{k+1} - y^*), \lambda^k - \lambda^*\rangle. 
\end{align}
Left-multiplying both sides of \eqref{yk+1-y*} by $\rho^{\frac{1}{2}} W^{\frac{1}{2}}$ yields, we have
\vspace{-0.3cm}
\begin{equation*}
\begin{aligned}
    & \rho^{\frac{1}{2}} W^{\frac{1}{2}} \left[(y^{k+1} - y^*) - \frac{1}{\eta} (\lambda^k - \lambda^*)\right] \\
    =&\rho^{\frac{1}{2}} W^{\frac{1}{2}} \left[\left(I_{nm}\!-\!\frac{\rho W}{\eta}\right)(y^k - y^*) + \frac{1}{\eta} A (x^{k+1} - x^*)\right].
\end{aligned}
\end{equation*}
Then, taking square both sides of the above equation results in
\vspace{-0.3cm}
\begin{align}\label{based_on_yupdate}
    & \rho \| y^{k+1} \!-\! y^* \|^2_W \!-\! 2 \left\langle \rho W ( y^{k+1} \!-\! y^*), \frac{1}{\eta} (\lambda^k \!-\! \lambda^*) \right\rangle \notag  \\
    & \!+\! \frac{\rho}{\eta^2} \| \lambda^k \!-\! \lambda^*\|^2_W  \notag \\
    =&\rho  \| y^k \!-\! y^*\|_{W\left(I_{nm}\!-\!\frac{\rho W}{\eta}\right)^2}^2 \!+\!\frac{\rho}{\eta} \| A(x^{k+1} \!-\! x^*)\|_W^2 \notag \\
    & \!+\!2 \left\langle \rho W\left(I_{nm}\!-\!\frac{\rho W}{\eta}\right)(y^k \!-\! y^*), \frac{1}{\eta} A(x^{k+1} \!-\! x^*)\right\rangle \notag \\
    =&\rho \| (y^k \!-\! y^*)\|_{W\left(I_{nm}\!-\!\frac{\rho W}{\eta}\right)^2}^2 \!+\!\frac{\rho}{\eta} \| A (x^{k+1} \!-\! x^*)\|_W^2 \notag \\
    & \!+\! \frac{2}{\eta}\left\langle \rho W\left(I_{nm}\!-\!\frac{\rho W}{\eta}\right)(y^k \!-\! y^*),  A(x^k \!-\! x^*)\right\rangle  \notag \\
    & \!-\!\frac{2}{\eta} \!\left\langle \! \! \rho W \! \!\left(\! I_{nm}\!-\!\frac{\rho W}{\eta} \!\right)\! \! (y^k \!-\! y^*), \! \alpha A (\!\nabla f(x^k) \!-\! \nabla f(x^*)\! ) \! \!\right\rangle  \notag \\
    & \!-\!\frac{2}{\eta} \left\langle \rho W\left(I_{nm}\!-\!\frac{\rho W}{\eta}\right)(y^k \!-\! y^*),  \alpha A A^{\rm{T}} ( y^k \!-\! y^*)\right\rangle \notag \\
    \leq& \rho  \| y^k \!-\! y^*\|_{W\left(I_{nm}\!-\!\frac{\rho W}{\eta}\right)^2}^2 \!+\!\frac{\rho}{\eta} \| A(x^{k+1} \!-\! x^*)\|_W^2  \notag \\
    & \!+\!\frac{1}{\eta} \| y^k \!-\! y^*\|_{(\theta_1 \!+\! \theta_2)\rho^2 W^2\left(I_{nm}\!-\!\frac{\rho W}{\eta}\right)^2 \!-\! \alpha \rho A A^{\rm{T}} W\left(I_{nm}\!-\!\frac{\rho W}{\eta}\right)}^2  \notag \\
    & \!+\! \frac{\| A \|^2}{\eta\theta_1} \| x^k \!-\! x^*\|^2 \!+\! \frac{\alpha^2\| A \|^2}{\eta\theta_2} \| \nabla f(x^k) \!-\! \nabla f(x^*) \|^2 \notag \\
    \leq & \| y^k \!-\! y^*\|_{H}^2 \!+\! \frac{\rho \| A \|^2 \! \| W \|}{\eta}  \!\| x^{k+1} \!-\! x^*\|_{}^2 \!+\! \frac{\| A \|^2}{\eta \theta_1} \| x^k \!-\! x^*\|^2 \notag \\
    & \!+\! \frac{\alpha^2 \| A \|^2}{\eta \theta_2}\| \nabla f(x^k) \!-\! \nabla f(x^*) \|^2,
\end{align}
where $\theta_1,~ \theta_2 > 0$ and $ H = \rho W \left(I_{nm}\!-\!\frac{\rho W}{\eta}\right)^2 \!+\! \frac{1}{\eta}\Big[(\theta_1 \!+\! \theta_2)\rho^2 W^2 \left(I_{nm}\!-\!\frac{\rho W}{\eta}\right)^2 \!-\! \alpha A A^{\rm{T}} \rho W  \left(I_{nm}\!-\!\frac{\rho W}{\eta}\right)\Big]$, the second equality results from \eqref{x update}, and the inequality is by the application of the triangular inequality. 

Combining \eqref{based_on_lambdaupdate} and \eqref{based_on_yupdate} leads to
\begin{align}\label{pk+1leqpk}
    &  \| y^{k+1} - y^* \|^2_{\rho W - \frac{\rho^2}{\eta} W^2} + \frac{1}{\eta} \| \lambda^{k+1} - \lambda^*\|^2 \notag \\
    & - \frac{\rho \| A \|^2 \| W \|}{\eta} \| x^{k+1} - x^*\|_{}^2 \notag \\
    \leq & \| \lambda^k - \lambda^*\|^2_{\frac{1}{\eta}I_{nm} - \frac{\rho}{\eta^2}W} + \| y^k - y^*\|_{H}^2 +\frac{\| A \|^2}{\eta \theta_1} \| x^k - x^*\|^2 \notag \\
    & +\frac{\alpha^2 \| A \|^2}{\eta \theta_2}\| \nabla f(x^k) - \nabla f(x^*) \|^2 .
\end{align}
Combining \eqref{omega_norm_1} and \eqref{pk+1leqpk} leads to 
\begin{align}\label{pk+1leq rhopk}
& \left(\!1 \!-\! \frac{\rho \| A \|^2 \!\| W \|}{\eta} \!\right)\!\| x^{k+1} \!-\! x^* \|_{}^2 \! + \!  \| \lambda^{k+1} \! - \!\lambda^* \|_{\frac{\alpha}{\rho}W^{\dagger} \! + \!  \frac{1}{\eta}I_{nm}}^2   \notag \\
& + \| y^{k+1} - y^* \|_{(\alpha \eta I_{nm}  +   \rho W)(I_{nm} - \frac{\rho}{\eta} W)}^2 \displaybreak[0] \notag \\
\leq & \, \left(1\! + \! \frac{ \| A \|^2}{\eta \theta_1} - \frac{2 \mu l_f \alpha}{l_f\! + \! \mu} \right)\| x^k - x^* \|_{}^2 \displaybreak[0] \notag \\
& \! + \!  \| y^k \!-\! y^* \|_{\alpha \eta (I_{nm} \!-\! \frac{\rho}{\eta} W) \! + \!  H}^2 \! + \!  \| \lambda^k \!-\! \lambda^* \|_{\frac{\alpha}{\rho}W^{\dagger} + \frac{1}{\eta}I_{nm}\! - \frac{\rho}{\eta^2}W}^2 \displaybreak[0] \notag \\
& -\left\|\Delta h^{k+1}\right\|_{\Omega}^2\! + \! \alpha \theta_0\left\|\Delta x^{k+1}\right\|^2 \displaybreak[0] \notag \\
& - \! \left( \! \frac{2 \alpha}{l_f\! + \! \mu}\!-\! \frac{\alpha}{\theta_0} \! -\! \frac{\alpha^2 \! \| A \|^2}{\eta \theta_2} \! \right) \! \left\|\nabla\! f \!\left(x^k\right)\!-\! \nabla \!f\!\left(x^{*}\right) \right\|^2\!.
\end{align}
% which is
% \begin{equation}
% \begin{aligned}
% & \| h^{k+1} - p^* \|_{\Theta_1}^2 \\
% \leq & \| h^k - p^* \|_{\Theta_2}^2 - \left\|\Delta h^{k+1}\right\|_{\Omega}^2+\alpha \theta_0\left\|\Delta x^{k+1}\right\|^2 \\
% & -\left(\frac{2 \alpha}{l_f+\mu}-\frac{\alpha}{\theta_0} - \frac{\alpha^2}{\eta \theta_2}\right)\left\|\nabla f\left(x^k\right)-\nabla f\left(x^{*}\right)\right\|^2,
% \end{aligned}
% \end{equation}
To prove that there exists a $\delta \in (0,1)$ such that \eqref{theorem 2} holds, we first show that the following inequalities hold,
\begin{subequations}\label{5_inequalities}
\begin{alignat}{5}
    & -\frac{\rho \| A \|^2 \| W \|}{\eta}\succ \frac{ \| A \|^2 }{\eta \theta_1} -\frac{2 \mu l_f \alpha}{l_f+\mu} ,\label{5-1} \displaybreak[0] \\
    & (\alpha \eta I_{nm} \! + \! \rho W\!) \! \left(\!\!I_{nm}\! - \! \frac{\rho}{\eta} W \!\!\right) \! \!\succ\!  \alpha \eta \! \left(\! \! I_{nm} \! - \! \frac{\rho}{\eta} W \!  \right)\! + \!H \! , \label{5-2} \displaybreak[0]\\
    & \frac{\alpha}{\rho}W^{\dagger} + \frac{1}{\eta} I_{nm} \succ \frac{\alpha}{\rho}W^{\dagger} + \frac{1}{\eta}I_{nm} - \frac{\rho}{\eta^2}W, \label{5-3} \displaybreak[0] \\
    & 1 - \alpha \theta_0 \geq 0 ,\label{5-4} \displaybreak[0] \\
    & \frac{2 \alpha}{l_f+\mu}-\frac{\alpha}{\theta_0} - \frac{\alpha^2 \| A \|^2}{\eta \theta_2} \geq 0 .\label{5-5}
\end{alignat}
\end{subequations}
Inequality \eqref{5-1} holds by the condition of the theorem, that is, $\eta > \frac{l_f + \mu}{2 \mu l_f \alpha} \left(\rho \| A \|^2 \| W \| + \frac{ \| A \|^2}{\theta_1}\right).$ 
\eqref{5-2} is equivalent to
\begin{equation*}
\begin{aligned}
    & \frac{\rho}{\eta} W\left(I_{nm}\!-\!\frac{\rho W}{\eta}\right)\\
    &\! \left(\! \alpha A A^{\rm{T}}\! + \!\rho W \!-\! (\theta_1 \!+ \!\theta_2)\rho W\left(I_{nm}\!-\!\frac{\rho W}{\eta} \right)  \!\right) \!\succ \!O_{nm \times nm},
\end{aligned}
\end{equation*}
which is true if the following inequality holds
\begin{equation*}
\begin{aligned}
    \alpha A A^{\rm{T}} + \rho(1 -\theta_1 - \theta_2) W + \frac{\rho^2 \left(\theta_1 + \theta_2\right) W}{\eta} \!\succ\! O_{nm \times nm}.
\end{aligned}
\end{equation*}
Since $W$ and $A A^{\rm{T}}$ are positive semidefinite, $\theta_1 + \theta_2 < 1$ ensures the above inequality and hence \eqref{5-2}.
% % we have the eigen-decomposition of $W$ as $W = U \Lambda U^{-1}$, where $U$ is a non-singular matrix and $\Lambda = \left[\begin{array}{cccc}
% \lambda_{n} & 0 & 0 & 0 \\
% 0 & \cdots & 0 & 0 \\
% 0 & 0 & \lambda_{2} & 0 \\
% 0 & 0 & 0 &0
% \end{array}\right]$. 
% Hence, we can rewrite \eqref{eq5-2} as
% \begin{equation}\label{eq5-2}
% \begin{aligned}
%     & \alpha I - \rho W - (\theta_1 + \theta_2)\rho W(I_{nm}\!-\!\frac{\rho W}{\eta})\\
%     = &\alpha I - \rho U \Lambda U^{-1} - (\theta_1 + \theta_2)\rho U \Lambda U^{-1} U (I-\frac{\rho}{\eta}\Lambda )U^{-1}\\
%     =& U[\alpha I - \rho \Lambda - (\theta_1 + \theta_2)(\rho  \Lambda - \frac{\rho^2 }{\eta}\Lambda^2) ]U^{-1}
% \end{aligned}
% \end{equation}
% % \begin{equation}
% % \begin{aligned}
% %     \alpha U U^{-1} - \rho U \Lambda U^{-1} - (\theta_1 + \theta_2)\rho U \Lambda U^{-1} U (I-\frac{\rho}{\eta}\Lambda )U^{-1} \succ 0
% % \end{aligned}
% % \end{equation}
% Setting 
% $
% \theta_1 + \theta_2 < \min_{i=2,\ldots,n} \frac{\rho \lambda_i+\alpha}{\rho \lambda_i - \frac{\rho}{\eta}\lambda_i} 
% $
% makes sure
% \begin{equation}
% \begin{aligned}
%     \alpha I - \rho \Lambda - (\theta_1 + \theta_2)(\rho  \Lambda - \frac{\rho^2 }{\eta}\Lambda^2) \succ O_p
% \end{aligned}
% \end{equation}
% which guarantees the satisfaction of \eqref{5-2}.
Inequality \eqref{5-3} is true due to the positiveness of both $
\rho$ and $\eta$.
Setting $
\alpha < \frac{2}{l_f + \mu}
$, $\theta_0 \in (\frac{l_f + \mu}{2}, \frac{1}{\alpha}]$ and $\eta \geq \frac{\alpha\|A\|^2}{(\frac{2}{l_f + \mu}-\frac{1}{\theta_0}) \theta_2}$ ensures \eqref{5-4} and \eqref{5-5}.
% Therefore, we have 
% $$
% \eta > \max \{ \frac{\alpha}{(\frac{2}{l_f + \mu}-\frac{1}{\theta_0}) \theta_2}, \frac{l_f + \mu}{2 \mu l_f \alpha} (\rho \lambda_{\min}^+(W) + \frac{1}{\theta_1})\}.
% $$
% $$
% \alpha < \frac{2}{l_f + \mu}
% $$
% and 
% $$
% \frac{\rho}{\eta} < \frac{1}{\lambda_{\max}(W)(W)}.
% $$

Based on inequalities \eqref{5-4} and \eqref{5-5}, we can eliminate the nonpositive terms on the right-hand side of \eqref{pk+1leq rhopk}, leading to the following recursive inequality that relates the $(k+1)$th and $k$th iterations,
\begin{equation}\label{compact_pk+1leqpk}
\left\|h^{k+1}-h^*\right\|_{\Theta}^2 \leq \left\|h^k-h^*\right\|_{\Theta^{\prime}}^2,
\end{equation}
where $\Theta$ and $\Theta'$ are block-diagonal matrices defined as $\Theta = {\rm{diag}}\{ \Theta_1, \Theta_2, \Theta_3\}$ and $\Theta^{\prime} = {\rm{diag}}\{ \Theta_1, \Theta_2, \Theta_3\}$
with individual blocks given by $\Theta_1 = I_{nm} - \frac{\rho\| A \|^2 \|W\|}{\eta}$, $\Theta_2 = (\alpha \eta I_{nm} + \rho W)\left(I_{nm} \!-\! \frac{\rho}{\eta} W\right)$, $\Theta_3 = \frac{\alpha}{\rho}W^{\dagger} + \frac{1}{\eta}I_{nm}$, $\Theta_1^{\prime} = (1+\frac{\| A \|^2}{\eta \theta_1} - \frac{2 \mu l_f \alpha}{l_f+\mu})I_{nm}$, $\Theta_2^{\prime}=\alpha \eta \left(I_{nm} \!-\! \frac{\rho}{\eta} W\right) + H$, and $\Theta_3^{\prime}=\frac{\alpha}{\rho}W^{\dagger} + \frac{1}{\eta}I - \frac{\rho}{\eta^2}W$.
Furthermore, in view of the fact that inequalities \eqref{5-1}, \eqref{5-2}, and \eqref{5-3} hold by proper parameter selections, we have
$\Theta \succ \Theta^{\prime} \succeq O_p.$ 
Consequently, there exists a scalar $\delta \in (0,1)$ such that
$\delta \Theta \succeq \Theta^{\prime} $. Substituting this bound into \eqref{compact_pk+1leqpk}, we obtain the following inequality,
$$
\begin{aligned}
\left\|h^{k+1}-h^*\right\|_{\Theta}^2 & \leq \delta \left\|h^k-h^*\right\|_{\Theta}^2 + \left\|h^k-h^*\right\|_{\Theta^{\prime} - \delta \Theta} \\
& \leq \delta \left\|h^k-h^*\right\|_{\Theta}^2,
\end{aligned}
$$
which completes the proof. \hfill $\Box$
\end{pf}
% The result in Theorem \ref{the: linear convergence} implies that both the primal and dual variables converge linearly, providing a stronger guarantee than the sublinear rates commonly observed in first-order distributed optimization methods. 

\section{Numerical Experiments} \label{sec: Simulation}
To validate the effectiveness of the proposed algorithm and to compare its performance with existing methods, we conduct a series of numerical experiments on a standard benchmark system. Specifically, the proposed method is compared with five state-of-the-art distributed algorithms: the Augmented Lagrangian Tracking (ALT) algorithm from \cite{falsone2023augmented}, Mirror-PG-EXTRA and Mirror-P-EXTRA from \cite{nedic2018improved}, DMAC from \cite{wu2025distributed}, and DuSPA from \cite{xu2018dual}. Notably, the proposed algorithm and Mirror-PG-EXTRA rely solely on gradient information, while Mirror-P-EXTRA, DMAC, and ALT require solving local optimization subproblems via the $\operatorname{argmin}$ operator at each iteration.

We consider the economic dispatch problem over the IEEE 118-bus system \cite{venkatesh2003comparison}, where the communication network is modeled as an undirected graph $\mathcal{G} = (\mathcal{V}, \mathcal{E})$ with $\mathcal{V} = {1, 2, \ldots, 118}$ and edges defined by $(i, i+1)$ and $(i, i+2)$ for $1 \leq i \leq 116$. This yields a sparse yet connected topology. A doubly stochastic weight matrix is generated accordingly. Among the 118 buses, 14 are randomly designated as generator buses, denoted by the set $\mathcal{V}_{\rm{g}}$. 
As for the settings for the optimization, we adapt the emission cost function in \cite{he2019optimizing} to model the following optimization problem, whose local cost of each generator is a combination of quadratic functions and exponential functions 
\begin{align}
    \min _{P} &~ f(P)= \sum_{i=1}^{118} a_i P_i^2+b_i P_i+\delta_i {\rm{e}}^{\ell_i P_i},  \\
    \operatorname{s.t.} &~  \sum_{i=1}^{118} P_i = d, \quad \underline{p}_i \leq P \leq \Bar{p}_i, \quad \forall i \in \mathcal{V}, \notag
\end{align}
where $P = [P_1\, P_2\, \ldots\, P_{118}]^{\rm{T}}$ and the coefficients for each generator are randomly generated positive numbers. To be more specific, for $i \in \mathcal{V}_{\rm{g}}$, $a_i \sim\mathcal{U}([0.3,0.7]$, $b_i\sim\mathcal{U}([100,400])$, $\delta_i \sim \mathcal{U}([1,10]) \times 10^{-4}$ and $l_i \sim \mathcal{U}([1,10]) \times 10^{-2}$. As for the limitations of generators, we set $[\underline{p}_i, \overline{p}_i] = [0,250]$, $i \in \mathcal{V}_{\rm{g}}$. For buses without generators, i.e., for $i\in \mathcal{V} \setminus \mathcal{V}_{\rm{g}}$, their corresponding coefficients $a_i$, $b_i$ and $\delta_i$ are set to zero and their local constraints are set as $\underline{p}_i = \overline{p}_i = 0$. 
The total power demand across the system is given by $ \sum_{i=1}^{14} d_i =950 \mathrm{MW}$, where the summation is unknown to each bus. Without loss of generality, we set a virtual initial local demand $d_i$ as $(950/14) \mathrm{MW}$ at each bus with a generator and $0 \mathrm{MW}$ at all other buses. The optimal solution $P^*$ of this problem is obtained by solving the problem centrally using a state-of-the-art solver MOSEK in YALMIP \cite{Lofberg2004}, and is used as a reference for evaluating convergence and accuracy. For the $\operatorname{argmin}$ opreation in  Mirror-P-EXTRA, DMAC, and ALT, we also invoke the solver MOSEK in YALMIP to get the solutions of local optimization problems. Our experimental setup uses a Laptop with an Intel Core i7-13700HX CPU, and the MATLAB version is R2023b.
 
\begin{figure}[!htb]
    \centering        
    \includegraphics[width=0.38\textwidth]{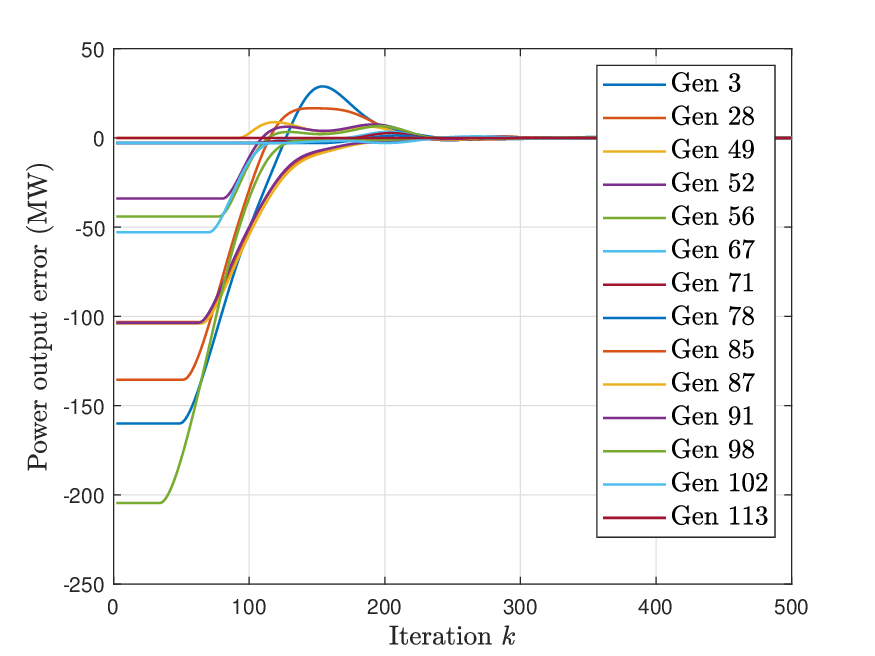}
    \caption{The power output error of each generator under the proposed Algorithm \ref{algorithm}.}
    \label{fig: power of each}
\end{figure}

\begin{figure}[!htb]
    \centering     
    \includegraphics[width=0.38\textwidth]{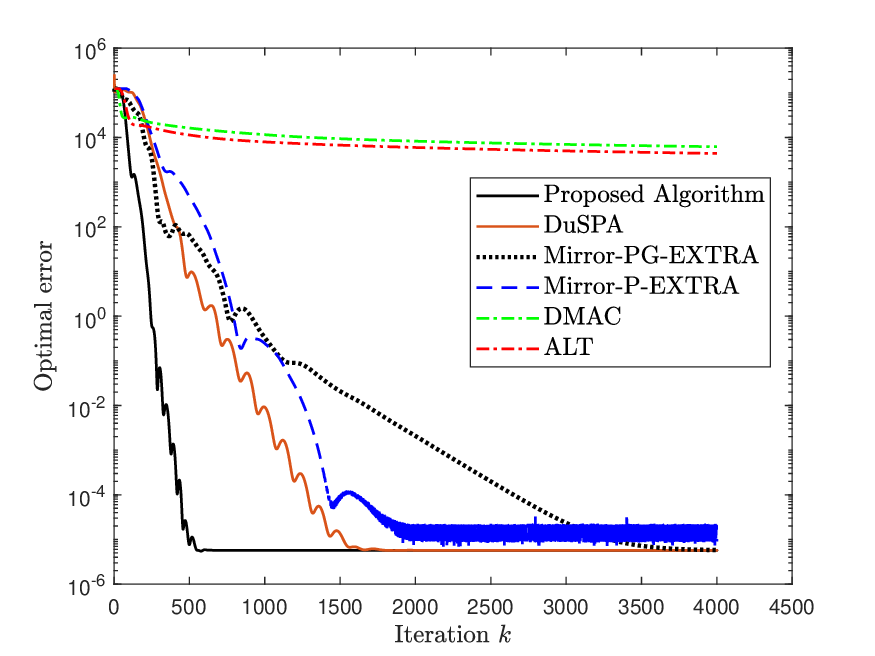}
    \caption{Comparison of optimal error over iteration for the six algorithms.}
    \label{fig: power mismatch}
\end{figure}

\begin{figure}[!htb]
    \centering     
    \includegraphics[width=0.38\textwidth]{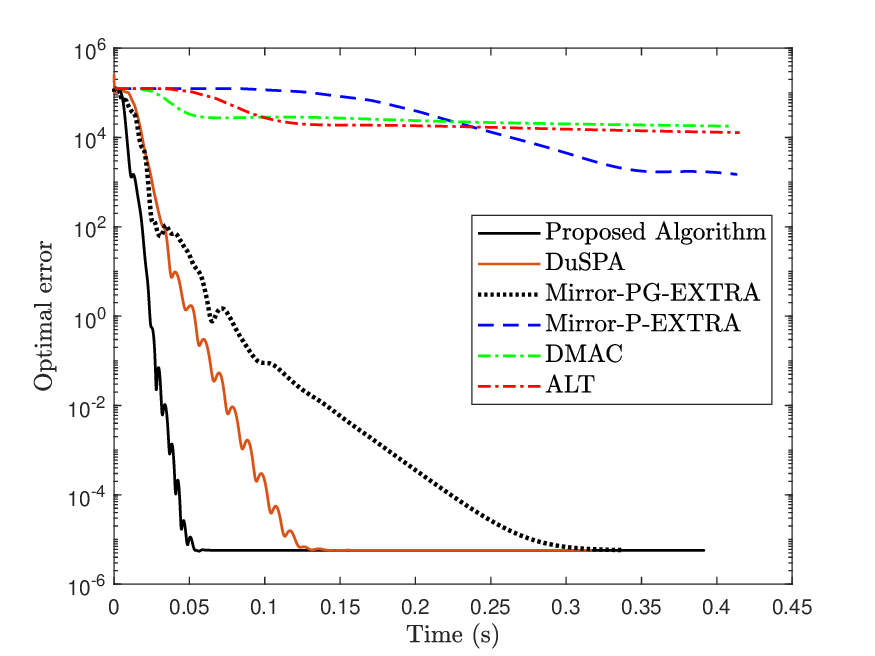}
    \caption{Comparison of optimal error over computation time for the six algorithms.}
    \label{fig: power mismatch time}
\end{figure}

% In relation to our algorithm, the variable $x$ denotes the power outputs of generators. When an optimal point is attained, $x$ gives the optimal outputs of generators.
Figure~\ref{fig: power of each} illustrates the evolution of the generator outputs, showing the difference between each $P_i$ and its optimal value $P_i^*$ over iterations. By the $300$th iteration, all generators exhibit negligible error, indicating convergence to the optimal dispatch.
Figure~\ref{fig: power mismatch} presents the power mismatch $\|P_k - P^*\|^2$ across iterations for the five algorithms. The proposed method consistently outperforms all baselines in terms of both convergence speed and solution quality. To assess computational time efficiency, Figure~\ref{fig: power mismatch time} plots the optimality error over actual computation time. It is observed that the proposed algorithm achieves the fastest convergence, significantly outperforming Mirror-P-EXTRA, Mirror-PG-EXTRA, DMAC, DuSPA, and ALT. This performance advantage arises because these methods require solving local optimization problems at each iteration. In contrast, the proposed algorithm, Mirror-PG-EXTRA, and DuSPA only involve first-order gradient computations and simple projections, resulting in superior time efficiency and faster error decay.

The simulation results on the IEEE 118-bus system demonstrate the effectiveness and efficiency of the proposed algorithm. Compared to five representative state-of-the-art methods, our algorithm achieves superior performance in terms of convergence speed and solution accuracy, both with respect to iteration count and CPU time. This advantage is primarily attributed to the algorithm’s lightweight per-iteration computation, which avoids solving local optimization subproblems and instead leverages gradient-based updates with simple projection steps. Furthermore, the proposed method maintains competitive accuracy while offering better scalability and practicality for resource-constrained distributed energy systems.
% \begin{figure}[!htb]
%     \centering        
%     \includegraphics[width=0.38\textwidth]{consensus_gradient.eps}
%     \caption{The consensus error for local gradients under the three algorithms.}
%     \label{fig: consensus error of gradient}
% \end{figure}
% \begin{figure}[!htb]
%     \centering        
%     \includegraphics[width=0.45\textwidth]{demand_violation.eps}
%     \caption{\xw{The absolute error of the total power generation and the total power demand under the five algorithms.}}
%     \label{fig: demand}
% \end{figure}

\vspace{-0.2cm}
\section{Conclusions}\label{sec: Conclusion}
\vspace{-0.2cm}
This paper presents a lightweight and scalable distributed algorithm for solving optimization problems with coupling equality constraints over networked systems. Unlike many existing methods that require solving local optimization problems at every iteration, the proposed algorithm obviates such a need by employing a projection-based gradient update. As a result, the proposed algorithm forms a single-loop iteration and reduces per-iteration computational overhead while preserving convergence guarantees.
% By reformulating the original problem as a dual consensus optimization problem, we develop an efficient first-order method that requires only local gradient evaluations and peer-to-peer communication.
We rigorously analyze the convergence behavior under both convex and strongly convex cost function settings, establishing nonergodic sublinear and linear convergence rates for general convex and strongly convex cases, respectively. Extensive numerical experiments on the IEEE 118-bus system validate the algorithm's effectiveness and demonstrate its superior performance in both convergence speed and computational efficiency compared to several benchmark methods. 
% These results confirm that the proposed method is well-suited for resource-constrained, and privacy-sensitive distributed scenarios.

\vspace{-0.2cm}
\bibliographystyle{unsrt}
\bibliography{reference}            
                                        % in the appendices.
\end{document}